\documentclass[11pt]{article}
\usepackage{amsmath}
\usepackage{amsfonts}
\usepackage{amsthm}
\usepackage{verbatim}

\begin{document}

\def\N{\mathbb{N}}
\def\F{\mathbb{F}}
\def\Z{\mathbb{Z}}
\def\R{\mathbb{R}}
\def\Q{\mathbb{Q}}
\def\H{\mathcal{H}}

\parindent= 3.em \parskip=5pt

\centerline{\bf{ ON THE CONTINUED FRACTION EXPANSION }}
\centerline{\bf{OF THE UNIQUE
    ROOT IN $\F(p)$ OF THE EQUATION }} 
\centerline{\bf{ $x^4+x^2-Tx-1/12=0$}}
\centerline{\bf{AND OTHER RELATED HYPERQUADRATIC EXPANSIONS}}

 \vskip 0.5 cm
\centerline{\bf{by A. Lasjaunias}}
\vskip 0.5 cm {\bf{Abstract.}}In 1985, Robbins observed by
computer the continued fraction expansion of certain algebraic power
series over a finite field. Incidentally he came across a particular
equation of degree four in characteristic $p=13$. This equation has
an analogue for all primes $p\geq 5$. There are two patterns for the continued fraction of the solution of this equation,
according to the residue of $p$ modulo 3. We describe this pattern
in the first case, considering especially $p=7$ and $p=13$. In the
second case we only give indications.   
\vskip 0.5 cm
\centerline{\bf{ 1. Introduction.}}
\par Throughout this note $p$ is an odd prime number, $\F_p$ is the finite field with $p$ elements and $\F(p)$ denotes 
the field of power series in $1/T$ with coefficients in $\F_p$, where $T$ is an indeterminate. These fields of power series are
known to be analogues of the field of real numbers. A non-zero
element of $\F(p)$ is represented by a power series expansion 
$$\alpha=\sum_{k\leq k_0}u_kT^k \quad \text{ where } k_0\in \Z, u_k\in
\F_p \quad \text{ and }u_{k_0}\neq 0.$$ 
We define $\vert \alpha\vert =\vert T\vert ^{k_0}$ where $\vert T\vert
>1$ is a fixed real number. The field $\F(p)$ is the completion of the
field $\F_p(T)$, of rational elements, for this absolute value.
\par Like in the case of real numbers, we recall that each irrational
element $\alpha \in\F(p)$, with $\vert \alpha\vert>1$, can be
expanded in an infinite continued fraction 
$$\alpha=[a_1,a_2,\dots,a_n,\dots]\quad \text{ where } a_i\in
\F_p[T]\text{ and }\deg(a_i)>0\text{ for }i\geq 1.$$  
The polynomials $a_i$ are called the partial quotients of the
expansion. For $n\geq 1$, we denote
$\alpha_{n+1}=[a_{n+1},a_{n+2},\dots]$, called the complete quotient, and we have
$$\alpha=[a_1,a_2,\dots,a_n,\alpha_{n+1}]=f_n(\alpha_{n+1})$$
where $f_n$ is a fractional linear transformation with coefficients in
$\F_p[T]$. Indeed, for $n\geq 1$, we have $f_n(x)=(x_nx+x_{n-1})/(y_nx+y_{n-1})$,
where the sequences of polynomials $(x_n)_{n\geq 0}$ and $(y_n)_{n\geq
  0}$, called the continuants, are both defined
by the same recursive relation: $K_n=a_nK_{n-1}+K_{n-2}$ for $n\geq 2$,
with the initials conditions $x_0=1$ and $x_1=a_1$ or $y_0=0$ and
$y_1=1$. Moreover, for $n\geq 1$, we have
$x_n/y_n=[a_1,\dots,a_n]$ and also $x_ny_{n-1}-x_{n-1}y_n=(-1)^n$.
\par We are
interested in describing the sequence of 
partial quotients for certain algebraic power series over
$\F_p(T)$. In the real case, an explicit description of
the sequence of partial quotients for algebraic numbers is only known
for quadratic elements. We will see that in the power series case over
a finite field, such a description is possible for many elements
belonging to a large class of algebraic power series containing the
quadratic ones. Our study is based upon a
particularly simple algebraic equation of degree 4.
\par Let $p$ be a prime number with $p\geq 5$. Let us consider the
following quartic equation with coefficients in $\F_p(T)$ :
$$x^4+x^2-Tx-1/12=0.\eqno{(1)}$$
 It is easy to see that $(1)$ has a unique root in $\F(p)$. We denote
 it by $u$ and we have $u=-1/(12T)+1/(12^2T^3)+\dots $. We put $\alpha=1/u$
and we consider the continued fraction expansion of $\alpha$ in
$\F(p)$. We have 
$$\alpha=[a_1,a_2,\dots,a_n,\dots] \quad \text{ with }\quad a_1=-12T.$$
A simple and general fact about this continued
fraction expansion can be observed. The root of $(1)$ is an odd function of $T$, since
$-u(-T)$ is also solution, and consequently all partial quotients are
odd polynomials in $\F_p[T]$. 
\par This quartic equation $(1)$ appeared for the first time in [6]. There the authors considered the case 
$p=13$, hence $-1/12=1$. A partial
conjecture on the continued fraction for the solution of (1) in $\F(13)$, observed by computer, 
was given in [6] and latter this conjecture was improved in [2]. The proof of this conjecture
 was given in [4]. 
\par The origin of the generalization for
arbitrary $p\geq 5$ replacing $1$ into $-1/12$ is to be found in [1]. If the continued fraction for the root of (1)
is peculiar and can be explicitely described, this is due to the
following result ([1],Theorem 3.1 , p. 263). 
\vskip 0.5 cm
\noindent {\emph{Let $p\geq 5$ be a prime number and
    $P(X)=X^4+X^2-TX-1/12 \in \F_p(T)[X]$. Then $P$ divides a
    nontrivial polynomial $AX^{r+1}+BX^r+CX+D$ where $(A,B,C,D)\in
    (\F_p[T])^4$ and where
$$r=p \quad \text{ if }\quad p=1 \mod 3 \quad \text{ and }\quad r=p^2 \quad \text{ if }\quad p=2 \mod 3.$$}} 
\par An irrational element $\alpha \in \F(p)$ is called
hyperquadratic if it satisfies  an algebraic
equation of the form $A\alpha^{r+1}+B\alpha^r+C\alpha+D=0$, where
$A,B,C$ and $D$ are in $\F_p[T]$ and $r$ is a power of $p$. The continued fraction expansion for many
hyperquadratic elements can be explicited. The reader may consult [8] for various examples and also more references.
\par According to the result
stated above, the solution in $\F(p)$ of $(1)$ is hyperquadratic. It appears
 that there are two different structures for the pattern of the continued
fraction expansion of this solution, corresponding to both cases : $p$
congruent to 1 or 2 modulo 3. In the second paragraph, if
$p=1 \mod 3$, we show that this continued fraction belongs to a much
larger family of hyperquadratic continued fractions and this allows us
to give an explicit description for $p=7$ and
$p=13$. In the second case we will only give indications which
might lead to an explicit description of the continued fraction. In a
last paragraph, we make a remark on programming which is based on a
result established by Mkaouar [7]. 
\par We need to introduce a pair of polynomials which play a fundamental
role in the expression of the continued fraction of the solution of
(1) and of many other algebraic power series. Throughout this note $k$ is an integer with $1\leq
k<p/2$. For $a \in \F_p^*$ we define a pair of
polynomials in
$\F_p[T]$ by 
$$P_{k,a}(T)=(T^2+a)^k\quad \text { and}\quad
Q_{k,a}(T)=\int_{0}^{T}P_{k-1,a}(x) dx.$$ 
Note that the definition of the second polynomial is made possible by
the condition $2k<p$ : by formal integration a primitive of $T^n$ in
$\F_p[T]$ is $T^{n+1}/(n+1)$ if $p$ does not divide $n+1$. 
\vskip 0.5 cm
\centerline{\bf{ 2. The case $p=1 \mod 3$.}}
\par Our method to describe the continued fraction expansion of the
solution of $(1)$, when $p=1 \mod 3$, is based upon the following conjecture.
\vskip 0.5 cm
\noindent {\bf{Conjecture 1.}}{\emph{Let $p$ be a prime number with $p=1 \mod
    3$. Let $\alpha\in \F(p)$ be defined by
    $-\alpha^4/12-T\alpha^3+\alpha^2+1=0$ and
    $\alpha=[a_1,a_2,\dots,a_n,\dots]$ its continued fraction
    expansion. Then there exist integers $l$ and $k$, a l-tuple $(\lambda_1,\lambda_2,\dots,\lambda_l)\in (\F_p^*)^l$ and a triple $(a,\epsilon_1,\epsilon_2) \in
(\F_p^*)^3$ such that 
$$(a_1,a_2,\dots,a_l)=(\lambda_1T,\lambda_2T,\dots,\lambda_lT)\eqno{(C_0)}$$ 
 and
$$\alpha^p=\epsilon_1P_{k,a}\alpha_{l+1}+\epsilon_2Q_{k,a}.\eqno{(C_1)}$$
We have $(l,k)=((p-1)/2,(p-1)/3)$.}}
\par We need to underline that the result stated above is only a conjecture
in the sense that it should be true for all primes $p$ with $p=1 \mod
3$. Actually for a particular prime $p$, a straightforward computation
implies $(C_0)$ and $(C_1)$. We will illustrate this for $p=7$ and $p=13$.
\par First we make the following observations. In [3] (Theorem 1 p. 332), it was proved that a unique power series in $\F(p)$ is well defined by
$(C_0)$ and $(C_1)$. Indeed, a unique power series in $\F(p)$ is well defined by
 $(C_1)$ and an arbitrary choice of the first $l$ partial quotients in
 $\F_p[T]$. Moreover this continued fraction is hyperquadratic. Indeed
  we have
    $$\alpha_{l+1}=(\alpha^p-\epsilon_2Q_{k,a})/\epsilon_1P_{k,a}\quad \text{
 and }\quad \alpha=(x_l\alpha_{l+1}+x_{l-1})/(y_l\alpha_{l+1}+y_{l-1})$$
    with the notations presented in the introduction. Combining these
    two equalities, we obtain the desired algebraic equation :
$$y_l\alpha^{p+1}-x_l\alpha^p+(\epsilon_1P_{k,a}y_{l-1}-\epsilon_2Q_{k,a}y_{l})\alpha+\epsilon_2Q_{k,a}x_{l}-\epsilon_1P_{k,a}x_{l-1}=0.$$
This conjecture is a stronger form of the theorem stated in the introduction.
\par {\emph{Proof of Conjecture 1 for $p=7$ and $p=13$:}} 
\newline Let
$\alpha$ be the inverse of the root of $(1)$. Then we have
$\alpha^4=-12(T\alpha^3-\alpha^2-1)$ and by iteration 
$$\alpha^n=a_n\alpha^3+b_n\alpha^2+c_n\alpha+d_n \quad \text{
  for }\quad n\geq 4,\eqno{(2)}$$
where $(a_n,b_n,c_n,d_n)\in (\F_p[T])^4$. It is easily checked that
for $p=7$ or $p=13$ we have $a_pb_{p+1}-a_{p+1}b_p=0$ in
$\F_p[T]$. Therefore combining $(2)$ for $n=p$ and for $n=p+1$, we
obtain
$$a_p\alpha^{p+1}-a_{p+1}\alpha^{p}=(a_pc_{p+1}-a_{p+1}c_p)\alpha + (a_pd_{p+1}-a_{p+1}d_p). \eqno{(3)}$$
We set $U_p=a_pd_{p+1}-a_{p+1}d_p$ and $V_p=a_{p+1}c_p-a_pc_{p+1}$,
consequently $(3)$ can be written as
$$\alpha=(a_{p+1}\alpha^{p}+U_p)/(a_{p}\alpha^{p}+V_p).\eqno{(4)}$$
We define $\delta$ as the g.c.d. of $a_p$ and $a_{p+1}$. Note that
$\delta$ is defined up to a multiplicative constant in $\F_p^*$. We
set $a_p^*=a_p/\delta$ and $a_{p+1}^*=a_{p+1}/\delta$. In the same way
, we set $U_p^*=a_p^*d_{p+1}-a_{p+1}^*d_p$ and
$V_p^*=a_{p+1}^*c_p-a_p^*c_{p+1}$. Consequently $(4)$ can be written as
$$\alpha=(a^*_{p+1}\alpha^{p}+U_p^*)/(a_{p}^*\alpha^{p}+V_p^*).\eqno{(5)}$$
We set $W=a_{p+1}^*V_p^*-a_p^*U_{p}^*$. For $p=7$ or $p=3$, we easily
check that $\vert a_{p}^*\alpha^p+V_p^*\vert > \vert
a_{p}^*W\vert$. Therefore, from $(5)$ we obtain
$$\vert \alpha -a^*_{p+1}/a_{p}^*\vert=\vert W\vert/\vert
a_{p}^*(a_{p}^*\alpha^p+V_p^*)\vert <\vert
a_{p}^*\vert ^{-2}.$$
This last inequality proves that $a^*_{p+1}/a_{p}^*$ is a convergent of
$\alpha$. For $p=7$ we have $a^*_{p+1}/a_{p}^*=[2T,6T,6T]=x_3/y_3$ and
for $p=13$ we have
$a^*_{p+1}/a_{p}^*=[T,12T,7T,11T,8T,5T]=x_6/y_6$. So the first part
of the conjecture holds for $p=7$ and $p=13$. As $\delta$ was
chosen up to a multiplicative constant, we can asume that
$a_{p+1}^*=x_l$ and $a_{p}^*=y_l$ with $l=(p-1)/2$ in both cases $p=7$
and $p=13$. Now we recall that we have
$$\alpha=(x_l\alpha_{l+1}+x_{l-1})/(y_l\alpha_{l+1}+y_{l-1}).\eqno{(6)}$$Hence, combining $(5)$ and $(6)$ we obtain 
$$\alpha^p=(-1)^lW\alpha_{l+1}+(-1)^l (x_{l-1}V_p^*-y_{l-1}U_{p}^*).\eqno{(7)}$$ 
By a simple computation , $(7)$ implies for $p=7$ 
$$\alpha^{7}=3(T^2-1)^2\alpha_4+5(5T^3+6T)\eqno
{(8)}$$
and for $p=13$
$$\alpha^{13}=(T^2+8)^4\alpha_7+4(2T^7+10T^5+12T^3+5T).\eqno {(9)}$$
So we see that the conjecture holds for $p=7$ with
$(\epsilon_1,\epsilon_2,a)=(3,5,6)$ and for $p=13$ with
$(\epsilon_1,\epsilon_2,a)=(1,4,8)$.  
\par We need to make a comment on the value of $a$ in the above conjecture
($6$ for p=7 and $8$ for $p=13$). When the paper [1] was prepared,
A. Bluher, interested in the Galois group of equation $(1)$,
could obtain some complementary results on the coefficients of the
hyperquadratic equation. At the fall of 2006, at a workshop in Banff,
she presented some of this work in progress. It results from these
formulas that we should have $a=8/27$ in all characteristic.
\par In order to normalize and to reduce the number of parameters, we make the following transformation.
 We define in
$\F_p[T]$, for $1\leq k<p/2$ as above, the following pair of
polynomials : 
$$P_k(T)=(T^2-1)^k\quad \text { and}\quad Q_k(T)=\int_{0}^{T}P_{k-1}(x) dx.$$ 
Let $v$ be a square root of $-a$ in $\F_p$ or $\F_{p^2}$, i.e. $v^2=-a$. Then we have
$P_{k,a}(vT)=(-a)^kP_k(T)$ and $Q_{k,a}(vT)=(-a)^{k-1}vQ_k(T)$. We put 
$\beta(T)=v\alpha(vT)$ and $\beta=[b_1,b_2,\dots,b_l,\beta_{l+1}]$. So if $\alpha=[a_1,a_2,\dots,a_l,\alpha_{l+1}]$ we obtain 
$\beta(T)=v\alpha(vT)=[va_1(vT),v^{-1}a_2(vT),\dots,v^{(-1)^l}\alpha_{l+1}(vT)]$. Therefore we have 
$\beta_{l+1}(T)= v^{(-1)^l}\alpha_{l+1}(vT)$ and
$b_i(T)=v^{(-1)^{i+1}}a_{i}(vT)$ for $i\geq 1$. Consequently $(C_1)$ can be written as 
$\alpha^p(vT)=(-a)^k\epsilon_1P_k\alpha_{l+1}(vT)+(-a)^{k-1}v\epsilon_2Q_k$. Finally $(C_1)$ becomes
$$\beta^p=\epsilon'_1P_k\beta_{l+1}+\epsilon'_2Q_k \eqno{(C'_1)}$$
where $$\epsilon'_1=(-a)^{k+(p-(-1)^l))/2}\epsilon_1 \quad \text{and}\quad  \epsilon'_2=(-a)^{k+(p-1)/2}\epsilon_2.$$
 While $(C_0)$ becomes 
$$(b_1,b_2,b_3,\dots,b_l)=(-a\lambda_1T,\lambda_2T,-a\lambda_3T,\dots).\eqno{(C'_0)}$$
To illustrate this transformation, which will be used later on, we apply it to the inverse of the solution of (1) in $\F(13)$.
Here we put $\beta=v\alpha(vT)$ with $v\in \F_{169}$ is such that
$v^2=5$. Consequently the first six partial quotients become
$$(b_1,b_2,b_3,b_4,b_5,b_6)=(5T,12T,9T,11T,T,5T) $$
 and (9) becomes 
$$\beta^{13}=12P_{4}\beta_7+9Q_{4}.\eqno
{(10)}$$
\par To describe the continued fraction expansion for the solution of (1) for $p=1 \mod 3$
, we need to introduce  a sequence 
of polynomials in $\F_p[T]$ based on the polynomial $P_k$. For a fixed
$k$ with $1\leq k<p/2$, we set
$$A_{0,k}=T \quad A_{i+1,k}=[A_{i,k}^p/P_k] \quad \text{ for }\quad 
i\geq 0.$$ 
Here the brackets denote the integral (i.e. polynomial) part of the rational function. 
We observe that the polynomials $A_{i,k}$ are odd polynomials in $\F_p[T]$. 
Moreover, it is important to notice that in the extremal case, if
$k=(p-1)/2$ then $A_{i,k}=T$ for $i\geq 0$.  
\par Now we shall consider all the continued fraction expansions defined by 
$$a_j=\lambda_jA_{i(j),k}  \quad \text{ for }\quad 1\leq j\leq l, \quad i(j) \in \N\quad \text{ and }\quad \lambda_j \in \F_p^* 
\leqno{(I)}$$ 
together with 
$$\alpha^p=\epsilon_1P_{k}\alpha_{l+1}+\epsilon_2Q_{k}\quad \text{
  where }\quad (\epsilon_1,\epsilon_2)\in (\F_p^*)^2.\leqno{(II)}$$ 
Here $p$ is an arbitrary odd prime, $l$ and $k$ are integers
 with $l\geq 1$ and $1\leq k<p/2$, and the plolynomials $P_k$, $Q_k$ and $A_{i,k}$ are defined in $\F_p[T]$ as above. We say that such an expansion is of type 
$(p,l,k)$. Note that the inverse of
the roots of (1) for $p=7$ or $p=13$, and conjecturally for all $p$
with $p=1 \mod 3$, have an expansion of type $(p,(p-1)/2,(p-1)/3)$ with
$a_j=\lambda_jA_{0,k}$ for $1\leq j\leq (p-1)/2$. Note that all power
series defined by a continued fraction of type $(p,l,k)$ are
hyperquadratic and they satisfy an algebraic equation of degree $p+1$. In the case of the root of $(1)$ this equation is reducible.  
\par The structure of expansions of type $(p,l,k)$ is based upon certain
properties of the pair $(P_{k},Q_{k})$ which are given in the
following proposition, the proof of which is to be found in [3]. 
\vskip 0.2 cm
\noindent {\bf{Proposition 1.}}{\emph{Let $p$ be an odd prime and $k$ an
  integer with $1\leq k<p/2$. We have in $\F_p(T)$ the following
  continued fraction expansion :
 $$P_k/Q_k=[v_{1,k}T,\dots,v_{i,k}T,\dots,v_{2k,k}T],$$
where the numbers $v_{i,k}\in \F_p^*$ are defined by $v_{1,k}=2k-1$ and recursively, for $1\leq i\leq 2k-1$, by
$$v_{i+1,k}v_{i,k}=(2k-2i-1)(2k-2i+1)(i(2k-i))^{-1}.$$
We set $\theta_k=(-1)^k\prod_{1\leq j\leq k}(1-1/2j)$.
Then we also have 
$$P_k/Q_k=-4k^2\theta_k^2[v_{2k,k}T,\dots,v_{1,k}T]\quad \text{ and
}\quad A_{i,k}^p=A_{i+1,k}P_k-2k\theta_k^{i+1}Q_k$$
for $i\geq 0$.}}
\par For certain expansions of type
 $(p,l,k)$, we have observed that all the partial
 quotients, obtained by computer, belong to the sequence $(A_{i,k})_{i\geq 1}$ up to
 a multiplicative constant in $\F_p^*$. In general the
 partial quotients are proportional to $A_{i,k}$ only up to a certain rank
depending on the choice of the first $l$
 partial quotients and of the pair $(\epsilon_1,\epsilon_2)$. Our goal was to 
understand under which conditions an expansion of type $(p,l,k)$ could
satisfy $a_j=\lambda_jA_{i(j),k}$ where $\lambda_j \in \F_p^*$ for all 
$j\geq 1$. In this case, we shall say that such an expansion is
perfect. We observe that if $k=(p-1)/2$, in a perfect expansion all
the partial quotients are proportional to $T$ : amazingly such an example also
exists in [6] (see the introduction of [3]).   
\par The following theorem gives a sufficient condition for an expansion of type $(p,l,k)$ to be perfect. Before stating our theorem, we need to describe the sequence $(i(n))_{n\geq 1}$ in $\N^*$ when the expansion is perfect.
 Given $l\geq1$ and $k\geq 1$, we define the sequence of integers $(f(n))_{n\geq 1}$ where
  $f(n)=(2k+1)n+l-2k$. Then the sequence of integers $(i(n))_{n\geq 1}$ is defined in the
following way :
$$i(j)\quad \text{is given for}\quad 1\leq j\leq l \quad \text{and }\quad i(f(n))=i(n)+1\quad \text{for }\quad n\geq 1,$$
$$i(n)=0 \quad \text{if }\quad n\notin f(\N^*)\quad \text{and }\quad
n>l.$$
\vskip 0.2 cm
\noindent {\bf{Theorem 1.}}{\emph{ Let $p$ be an odd prime, $k$ and
    $l$ be given as above. Let
    $(\lambda_1,\dots,\lambda_l)$ in $(\F_p^*)^l$ and
    $(i(1),\dots,i(l))$ in $\N^l$ be given. Let $(f(n))_{n\geq 1}$ and
    $(i(n))_{n\geq 1}$ sequences of integer, $(A_{i,k})_{i\geq 0}$
    sequence in $\F_p[T]$ and $\theta_k\in \F_p^*$ be defined as above. Let $\alpha \in \F(p)$ be a continued fraction of type $(p,l,k)$
 defined by $(I)$ and $(II)$. 
 If we can define in $\F_p^*$
$$(III)\qquad  \delta_n=[\theta_k^{i(n)}\lambda_n,\dots,\theta_k^{i(1)}\lambda_1,2k\theta_k/\epsilon_2] \quad
\text{ for }\quad 1\leq n\leq l$$ 
and we have 
$$(IV)\qquad\delta_l=2k\epsilon_1/\epsilon_2$$
then the partial quotients of this expansion satisfy 
$$(V)\quad a_n=\lambda_nA_{i(n),k}\quad \text{ where }\quad \lambda_n \in
\F_p^* \quad \text{ for }\quad n\geq 1.$$
Moreover the sequence $(\lambda_n)_{n\geq 1}$ in
$\F_p^*$ is defined recursively by the first values $\lambda_1,\lambda_2,\dots,\lambda_l$ and 
 for $n\geq 1$ : $\lambda_{f(n)}=\epsilon_1^{(-1)^{n}}\lambda_n$, 
$$\lambda_{f(n)+i}
=-v_{i,k}\epsilon_1^{(-1)^{n+i}}(2k\theta_k\delta_n)^{(-1)^i} \quad \text{ for }1\leq i\leq 2k$$
together with $(\delta_n)_{n\geq 1}$ defined recursively in $\F_q^*$ by the
initial values $\delta_1,\dots,\delta_l$ given in $(III)$ and for
$n\geq 1$:  $\delta_{f(n)}=\epsilon_1^{(-1)^{n}}\delta_n\theta_k$, 
 $$\delta_{f(n)+i}=\epsilon_1^{(-1)^{n+i}}(iv_{i,k}/(2k-2i+1))(2k\theta_k\delta_n)^{(-1)^i}
 \quad \text{ for }1\leq i\leq 2k.$$
}}
\par This theorem is a modified version of a stronger one given in [5] (Theorem B, p.256). The first modification
is due to a simplification. In our previous works, we have considered a more general situation where the power series
 are defined over a finite field $\F_q$ not necessarily prime. Also we introduced a larger class of continued fraction expansions
of type $(r,l,k)$ where $r$ is a power of $p$, and $k$ is chosen in a particular subset of integers in relation with $r$.
In this more general context, condition $(IV)$ of Theorem 1 is only sufficient to have all partial quotients proportional to $A_{i,k}$
(see [5] Corollary C). However, if the base field is prime, we may think that $(III)$ and $(IV)$ are sufficient and necessary conditions
 to have a perfect expansion. Now we need to explain the main modification. Earlier we had assumed that the $l$ first partial quotients 
were proportional to $A_{0,k}=T$. Here we make a larger hypothesis, which is $(I)$. Indeed this new hypothesis,
 does not alter the proof of Theorem B to which the reader is referred. Nevertheless this 
implies a minor change: condition $(II)$ there in Theorem B ([5], p. 256) becomes $(III)$ here in Theorem 1. 
\par We illustrate Theorem 1 in the simple case where $l=k=1$ :
\vskip 0.2 cm
\noindent {\bf{Corollary 1.}}{\emph{Let $p$ be an odd prime. Let
    $i(1)\in \N$ and $(i(n))_{n\geq 1}$ be defined recursively in $\N$ by 
$$i(3n-1)=i(n)+1 \quad \text{and}\quad i(3n)=i(3n+1)=0\quad
\text{for}\quad  n\geq 1.$$
Let $(A_{i,1})_{i\geq 1}$ be the sequence in $\F_p[T]$ defined above.
Let $p$ be an odd prime and $\alpha \in \F(p)$ a continued fraction
 of type $(p,1,1)$, with $a_1=\lambda_1A_{i(1),1}$. Assume that $\epsilon_2^2+2\epsilon_1\neq 0$
and $\lambda_1=(\epsilon_2^2+2\epsilon_1)(-2)^{i(1)}/\epsilon_2$. 
Then we have $a_n=\lambda_nA_{i(n),1}$ where the sequence $(\lambda_n)_{n\geq 2}$ in $\F_p^*$ is defined by 
$$\lambda_{3n-1}=\epsilon_1^{(-1)^n} \lambda_n, \quad \lambda_{3n}=-\epsilon_1^{(-1)^{n+1}} \delta_n^{-1} ,\quad \lambda_{3n+1}=-\lambda_{3n}^{-1}$$
and  the sequence $(\delta_n)_{n\geq 1}$ in $\F_p^*$ is defined by
$\delta_1=-2\epsilon_1/\epsilon_2$ and 
$$\delta_{3n-1}=-\epsilon_1^{(-1)^n} \delta_n/2,\quad
\delta_{3n}=-\epsilon_1^{(-1)^{n+1}} \delta_n^{-1}, \quad
\delta_{3n+1}=2\delta_{3n}^{-1}.$$}}
\par Now we turn our attention to the continued fraction expansion of
the root of equation (1). By applying Theorem
1, we have the following result.
\vskip 0.2 cm 
\noindent {\bf{Corollary 2.}}{\emph{Let $p$ be an odd prime with $p=1
    \mod 3$.  We set $(l,k)=((p-1)/2,(p-1)/3)$. We consider the
    sequence $(A_{i,k})_{i\geq 0}$ in $\F_p[T]$ introduced above. For $n\geq 1$, we set $i(n)=v_{(2p+1)/3}((p-1)(4n-1)/6)$ where $v_m(n)$ denotes the largest
power of $m$ dividing $n$. Let $\alpha$ be the unique root of $-x^4/12-Tx^3+x^2+1=0$
    in $\F(p)$ and $\alpha=[a_1,a_2,\dots,a_n,\dots]$ its continued
    fraction expansion. Assume $p=7$ or $p=13$ and set $v=1$ if $p=7$ and $v=\sqrt{5}\in
\F_{169}$ if $p=13$. Then there exists a sequence
    $(\lambda_n)_{n\geq 1}$ in $\F_p^*$ such that
$$a_n=\lambda_nv^{(-1)^n}A_{i(n),k}(T/v)\quad \text{ for } n\geq 1.$$  
 The sequence $(\lambda_n)_{n\geq 1}$ is
defined from the $(l+2)$-tuple
$(\epsilon_1,\epsilon_2,\lambda_1,\dots,\lambda_l)$, as in Theorem 1, with
$(3,5,2,6,6)$ for $p=7$ and $(12,9,5,12,9,11,1,5)$ for $p=13$. 
}}
\par The proof of this corollary follows immediately from the proof of
Conjecture 1 for $p=7$ or $p=13$. We put $\beta(T)=v\alpha(vT)$. Hence
we have $\beta=[b_1,\dots,b_n,\dots]$ where
$a_n(T)=v^{(-1)^{n}}b_n(T/v)$. If
$p=7$ then the continued fraction expansion for $\beta$ is defined by 
$(b_1,b_2,b_3)=(2T,6T,6T)$ and $(8)$, moreover $k=2$ and $\theta_2=3$. If
$p=13$ then the continued fraction expansion for $\beta$ is defined by 
$(b_1,b_2,b_3,b_4,b_5,b_6))=(5T,12T,9T,11T,T,5T)$ and $(10)$, moreover
$k=4$ and $\theta_4=2$. To apply Theorem 1, we only need to check that
$[\lambda_l,\lambda_{l-1},\dots,\lambda_1,2k\theta_k/\epsilon_2]$
exists and is equal to $2k\epsilon_1/\epsilon_2$, in both cases. Hence we
have $b_n=\lambda_nA_{i(n),k}$ for $n\geq 1$. Finally it is elementary
to verify that the sequence $(v_{(2p+1)/3}((p-1)(4n-1)/6))_{n\geq 1}$
satisfies the same initial conditions and the same recurrence relation
as the sequence $(i(n))_{n\geq 1}$ defined before Theorem 1.
\par As we remarked after Conjecture 1, the limitation to $p=7$ and $p=13$ in this corollary is
artificial. Indeed, by computer, for a given $p=1 \mod 3$, after the
transformation mentioned above where $v=\sqrt{-8/27}$, it is
possible to obtain the $(l+2)$-tuple,
$(\epsilon_1,\epsilon_2,\lambda_1,\dots,\lambda_l)$. Thereby we can check that
the right condition is fullfilled. Therefore we conjecture that the
formula given in this corollary for the partial quotients of the root
of $(1)$ holds for all primes $p$ with $p=1 \mod 3$. Besides, we recall that the case $p=13$ of this corollary has already been
published in [4].
\par  To measure the quality of rational approximation to a given
irrational power series, we have the following classical definition. 
\vskip 0.2 cm  
\noindent {\bf{Definition.}}{\emph{Let $\alpha \in \F(p)$ be an irrational element and $\alpha=[a_0,a_1,\dots,a_n,\dots]$ its continued 
fraction expansion. We set $$\nu_0(\alpha)=\limsup_{n\geq 0}(\deg(a_{n+1})/\sum_{0\leq k\leq n}\deg(a_k))\in \R^+ \cup \lbrace \infty \rbrace .$$
The quantity $\nu(\alpha)=2+\nu_0(\alpha)$ is called the rational approximation exponent of $\alpha$. }}
\par Note that an analogous quantity for real numbers can be
defined. In the middle of the nineteenth century, Liouville remarked
that this quantity was bounded for algebraic numbers and so he could prove the existence of transcendental real numbers.
 In the middle of the twentieth century, Mahler adapted Liouville's work to the
setting of fields of power series over an arbitrary field. If $\alpha$ is algebraic of degree $d$ over
 $\F_p(T)$ then we have $\nu(\alpha)\in [2;d]$. If $\alpha \in \F(p)$
 is defined as a continued fraction expansion of type $(p,l,k)$, we know that
it satisfies an algebraic equation of degree $p+1$ and consequently we have $\nu(\alpha)\in [2;p+1]$. If this expansion 
is perfect, the description given in Theorem 1 allows to compute the rational approximation
exponent. Indeed, since we have $\deg
A_{0,k}=1$ and $\deg(A_{i+1,k})=p\deg(A_{i,k})-2k$, we see that
$a_n=\lambda_nA_{i(n),k}$ implies
$\deg(a_n)=(p^{i(n)}(p-1-2k)+2k)/(p-1)$. Consequently, if the sequence $i(n)$ is
not too complex, the computation of $\nu_0(\alpha)$ is elementary. We
just state below two cases.
\vskip 0.2 cm  
\noindent {\bf{Corollary 3.}}{\emph{Let $\alpha \in \F(p)$ be a perfect continued fraction of type $(p,l,k)$.
\newline If  $a_j=\lambda_jA_{0,k}$ for $1\leq j\leq l$ then $\nu_0(\alpha)=(p-2k-1)/l$.
\newline If  $l=k=1$ and $a_1=\lambda_1A_{i,1}$ then $\nu_0(\alpha)=(p-1)(p^{i+1}-3p^i)/(p^{i+1}-3p^i+2)$.
Consequently we have $\nu(\alpha)=8/3$, if $\alpha$ is the root of
$(1)$ for $p=7$ or $p=13$.}}
\par Note that the second result in this corollary implies that a perfect continued fraction
of type $(p,1,1)$ is algebraic of degree $p+1$ over
$\F_p(T)$, if $p\geq 5$ and $i\geq 1$. Moreover, according to the previous remarks, the last statement is conjectured
to be true for all $p=1 \mod 3$. 
\par Finally we make a remark on perfect expansions. For a given triple $(p,l,k)$ and a given vector $(i(1),\dots,i(l))\in \N^l$ 
there are clearly $(p-1)^{l+2}$ expansions of type $(p,l,k)$, each one defined by the $l+2$-tuple 
$(\epsilon_1,\epsilon_2,\lambda_1,\dots,\lambda_l)$. A general 
expansion of type $(p,l,k)$ has a pattern difficult to describe and we have not tried to do so. Nevertheless, in such a 
finite set of expansions, it seems that the approximation 
exponent should be minimal when the expansion is perfect.
\vskip 0.5 cm
\centerline{\bf{ 3. Indications in the case $p=2 \mod 3$.}}
\par In this second case, the pattern for the continued fraction of the solution of equation
$(1)$ appears to be very different from the one we have described in
the first case. Here again there seems to be a general
pattern for all primes $p$ with $p=2 \mod 3$. This pattern is not
understood but we can indicate some observations which are somehow
parallel to what has been presented above. We do not know wether a similar method can be developped from these indications to
obtain an explicit description of this continued fraction. We have the
following conjecture.  
\vskip 0.2 cm
\noindent {\bf{Conjecture 2.}}{\emph{ Let $p$ be a prime number with
    $p=2 \mod 3$. Let $\alpha\in \F(p)$ be defined by
    $-\alpha^4/12-T\alpha^3+\alpha^2+1=0$ and
    $\alpha=[a_1,a_2,\dots,a_n,\dots]$ its continued fraction
    expansion. Then there exist integers $l$ and $k$ and a triple $(a,\epsilon_1,\epsilon_2) \in
(\F_p^*)^3$ such that 
$$\alpha^{p^2}=\epsilon_1P_{k',a}\alpha_{l+1}+\epsilon_2Q^p_{k,a}\eqno{(C_2)}$$
We have $$(l,k',k)=(\frac{(p+1)^2}{3},\frac{p^2-1}{3},\frac{p+1}{3}).$$}}
\par At last we state the following proposition, whose proof could be
simply deduced from Proposition 1. The notations are the same as there.
\vskip 0.2 cm
\noindent {\bf{Proposition 2.}}{\emph{Let $p$ be an odd prime, $k$ an
  integer with $1\leq k< p/2$ and $i$ an
  integer with $1\leq i< p/2$ . We have in $\F_p(T)$ the following
  continued fraction expansion :
 $$P_{kp-i}/Q_k^p=[v_{1,k}A_{1,i},-\delta_1^{-1}v_{1,i}T,
 -\delta_1v_{2,i}T\dots,-\delta_1^{-1}v_{2i-1,i}T,-\delta_1v_{2i,i}T,$$
$$\phantom{ P_{kp-i}/Q_k^p=[}v_{2,k}A_{1,i},-\delta_2^{-1}v_{1,i}T,
 -\delta_2v_{2,i}T\dots,-\delta_2^{-1}v_{2i-1,i}T,-\delta_2v_{2i,i}T,$$
$$\phantom{ P_{kp-i}/Q_k^p=[}\dots \quad \dots \quad \dots \quad \dots \quad \dots$$
$$,v_{2k-1,k}A_{1,i},-\delta_{2k-1}^{-1}v_{1,i}T,
 -\delta_{2k-1}v_{2,i}T\dots,-\delta_{2k-1}v_{2i,i}T,v_{2k,k}A_{1,i}]$$
where $$\delta_j=2i\theta_i[v_{j,k},v_{j-1,k},\dots,v_{1,k}]\quad \text{ for }\quad 
1\leq j\leq 2k-1.$$
Moreover, writting $[b_1,b_2,\dots,b_n]$ for the expansion given above, we
also have 
$[b_1,b_2,\dots,b_n]=-4k^2\theta_k^2[b_n,b_{n-1},\dots,b_1].$
}}
\vskip 0.5 cm
\centerline{\bf{ 4. A remark on programing.}}
\par Before concluding, we want to discuss a particular way to obtain
by computer the begining of the continued fraction expansion for an algebraic power
series. The natural way is to start from a rational approximation, often
obtained by tuncating the power series expansion, and therefrom
transform this rational into a finite continued fraction as this is done
for an algebraic real number. However, here in the formal case, it is
possible to process differently. The origin of this method is based on
a result introduced by M. Mkaouar, it can be found in [7] and also in other papers from
him. We recall here this result : 
\newline {\bf{Proposition }}(Mkaouar){\emph{Let $P$ be a polynomial in
    $\F_q[T][X]$ of degree $n\geq 1$ in $X$. We put $P(X)=\sum_{0\leq i\leq n}
    a_iX^i$ where $a_i\in \F_q[T]$ . Assume that we have 
$$\leqno{(*)}\qquad \vert a_i \vert < \vert a_{n-1} \vert \quad \text{ for }\quad 0\leq
i\leq n \quad \text{ and }\quad i\neq
    n-1.$$
      Then $P$ has a unique root in
    $\F(q)^+=\lbrace \alpha \in \F(q) \mid \vert \alpha \vert \geq \vert T \vert \rbrace$. Moreover, if $u$ is this root, we have
    $[u]=-[a_{n-1}/a_n]$. If $u\neq [u]$ and $u=[u]+1/v$ then
    $v$ is the unique root in
    $\F(q)^+$ of a polynomial $Q(X)=\sum_{0\leq i\leq n}
    b_iX^i$ with the same property $(*)$ on the coefficients $b_i$.}} 
\par In this proposition, it is clear that the coefficients $b_i$ can be
deduced from $[u]$ and the $a_i$'s. We have $b_n =P([u])$, so if $u$
is not integer we obtain $[v]=-[b_{n-1}/b_n]$. Consequently the process
can be carried on, for a finite number of steps if the solution $u$ is
rational or infinitely otherwise. Thus the partial quotients of the
solution can all be obtained by induction. This method can be applied
to obtain the continued fraction expansion of the solution of our
quartic equation, starting from the polynomial
$P(X)=-X^4/12-TX^3+X^2+1$. We have written here bebow the few lines of
a program (using Maple) to obtain the first two hundred partial quotients of this expansion.

\begin{verbatim}
p:=5:n:=200:u:=-1/12 mod p:
a:=array(1..n):b:=array(1..n):c:=array(1..n):d:=array(1..n):
e:=array(1..n):qp:=array(1..n):a[1]:=u:
b[1]:=-T:c[1]:=1:d[1]:=0:e[1]:=1:qp[1]:=-quo(b[1],a[1],T) mod p:
for i from 2 to n do 
a[i]:=simplify(a[i-1]*qp[i-1]^4+b[i-1]*qp[i-1]^3+
c[i-1]*qp[i-1]^2+d[i-1]*qp[i-1]+e[i-1]) mod p:
b[i]:=simplify(4*a[i-1]*qp[i-1]^3+3*b[i-1]*qp[i-1]^2+
2*c[i-1]*qp[i-1]+d[i-1]) mod p:
c[i]:=simplify(6*a[i-1]*qp[i-1]^2+3*b[i-1]*qp[i-1]+c[i-1]) mod p:
d[i]:=simplify(4*a[i-1]*qp[i-1]+b[i-1]) mod p:e[i]:=a[i-1]:
qp[i]:=-quo(b[i],a[i],T) mod p:od:print(qp);
\end{verbatim}

\noindent This method is easy to use in two cases : if the degree of
the initial polynomial is small and also if the initial polynomial has
the particular form corresponding to an hyperquadratic
solution. Indeed in this second case, in the proposition stated above,
the polynomial $Q$ has the same form as $P$, therefore the recurrence
relations between the coefficients of $P$ and those of $Q$ are made
simple. In both cases, it seems that the method is of limited
practical use because the degrees of the polynomials in $T$, giving the
partial quotients by division, are growing fast.

\noindent Lasjaunias Alain
\newline Institut de Mathématiques de Bordeaux-CNRS UMR 5251 
\newline Universit\'e Bordeaux 1 
\newline 351 Cours de la Lib\'eration
\newline F-33405 TALENCE Cedex  FRANCE
\newline e-mail: Alain.Lasjaunias@math.u-bordeaux1.fr


\begin{thebibliography}{7}    

\addcontentsline{toc}{section}{References}


\bibitem[1]{} A. Bluher and A. Lasjaunias,\emph{ Hyperquadratic power
    series of degree four}, Acta Arithmetica {\bf 124}
  (2006), 257-268.

\bibitem[2]{} W. Buck and D. Robbins,\emph{ The continued fraction of an algebraic power series 
satisfying a quartic equation}, Journal of Number Theory {\bf 50} (1995), 335--344.

\bibitem[3]{} A. Lasjaunias,\emph{ Continued fractions for
    hyperquadratic power series over a finite field}, Finite Fields
  and their Applications {\bf 14} (2008), 329-350. 

\bibitem[4]{} A. Lasjaunias,\emph{ On Robbins' example of a continued
    fraction expansion for a quartic power series over $\F_{13}$},
  Journal of Number Theory {\bf 128} (2008), 1109-1115.
 
\bibitem[5]{} A. Lasjaunias, \emph{ Algebraic continued fractions in
    $\F_q((T^{-1}))$ and recurrent sequences in $\F_q$ },  Acta Arithmetica {\bf 133}
  (2008), 251-265.  

\bibitem[6]{} W. Mills and D. Robbins,\emph{ Continued fractions for certain algebraic power 
series}, Journal of Number Theory {\bf 23} (1986), 388--404.

\bibitem[7]{} M. Mkaouar,\emph{ Sur les fractions continues des
    s\'eries formelles quadratiques sur $\F_q(X)$}, Acta Arithmetica {\bf 97.3}
  (2006), 241-251.

\bibitem[8]{} W. Schmidt,\emph{ On continued fractions and diophantine
    approximation in power series fields}, Acta Arithmetica {\bf 95.2}
  (2000), 139-166.

\end{thebibliography}
\end{document}